\input amstex
\documentstyle{amsppt}
\magnification=\magstep 1
\pageheight{21truecm}
\TagsOnRight
\CenteredTagsOnSplits

\topmatter
\title  Quasifinite representations of $W_{\infty}$ \endtitle
\rightheadtext{ Quasifinite representations of $W_{\infty}$ }
\leftheadtext{Victor G. Kac and Jos\'e I. Liberati}

\author{Victor G. Kac and Jos\'e I. Liberati} \endauthor 

\address  Department of Mathematics, MIT, Cambridge, MA 02139, USA \endaddress
\email kac\@math.mit.edu  \endemail
\address FAMAF, Universidad Nacional de C\'ordoba - (5000) C\'ordoba, 
Argentina.
\endaddress
\curraddr  Department of Mathematics, MIT, Cambridge, MA 02139, USA  
\endcurraddr 
\email liberati\@math.mit.edu  -  liberati\@mate.uncor.edu \endemail


  
\abstract  We classify the quasifinite highest weight modules over a family
of subalgebras  $W_{\infty}^{(n)}$ of the central extension $W_{1+\infty}$
of the Lie algebra of differential operators on the circle consisting of
operators of order  $\geq n$. We classify the unitary quasifinite highest
weight modules over  $W_{\infty}=W_{\infty}^{(1)}$ and  realize them  in
terms of unitary highest weight representations of  the Lie algebra of
infinite  matrices with  finitely many non-zero diagonals.

\endabstract
\endtopmatter
\document

\head 1. Introduction\endhead

The $\Cal W$-infinity algebras naturally arise in various physical
theories, such as conformal field theory, the theory of the quantum Hall
effect, etc.  The ${W}_{1+\infty}$ algebra, which  is the central extension
of the Lie algebra $\Cal D$ of differential operators  on the circle, is
the most fundamental among these algebras.

When we study the representation theory of a Lie algebra of this kind, we
encounter the difficulty that although it admits a $\Bbb Z$-gradation, each
of the graded subspaces is still infinite dimensional, and therefore  the
study of highest weight modules   which satisfy the quasi-finiteness
condition that its  graded subspaces  have finite dimension, becomes a
non-trivial problem.

 The study of   representations  
 of the Lie algebra ${W}_{1+\infty}$ was initiated in [KR1], where a
characterization of its  irreducible quasifinite highest weight
representations was given,  these modules were constructed in terms of
irreducible  highest weight representations of the Lie algebra of infinite
matrices, and  the unitary ones were described. On the basis of this
analysis, further studies were made in the framework of vertex algebra
theory for the ${W}_{1+\infty}$ algebra [FKRW, KR2], and for its matrix
version [BKLY]. The case of orthogonal subalgebras of $W_{1+\infty}$ was
studied in [KWY]. The symplectic subalgebra of ${W}_{1+\infty}$ was
considered in [B] in relation to number theory.

In this article, we develop the representation theory of the   subalgebras
$W_{\infty}^{(n)}$  of  $W_{1+\infty}$, where  $W_{\infty}^{(n)}$ ($n\in
\Bbb N$) is  the central extension of the Lie algebra $\Cal D \partial_t^n$
of differential operators on the circle that annihilate all polynomials of
degree less than $n$. The most important of these subalgebras is
$W_{\infty}=W_{\infty}^{(1)}$.  We first give a description of parabolic
subalgebras of $W_{\infty}^{(n)}$ and we classify  all its irreducible
quasifinite highest weight modules. Second,  we describe the relation of
$W_{\infty}$ to  the central extension of the Lie algebra of infinite
matrices with  finitely many non-zero diagonals and, using this relation,
we  classify and construct all the  unitary  irreducible quasifinite
modules over $W_{\infty}$. Surprisingly, the list of unitary modules over
$W_{\infty}$ is much richer than that over $W_{1+\infty}$.

\

%
%
%
%

\head 2. Quasifinite representations of $\Bbb Z$-graded Lie algebras \endhead

Let $\frak g$ be a $\Bbb Z$-graded Lie algebra over $\Bbb C$: 
$$
\frak g = \bigoplus_{j\in \Bbb Z} \frak g_j \ , \quad \quad [\frak g_i ,
\frak g_j ]\subset \frak g_{i+j},
$$
where $\frak g_i$ is not necesarily of finite dimension. Let $\frak g_{\pm}
= \oplus_{j>0} \frak g_{\pm j}$. A subalgebra $\frak p$ of $\frak g$ is
called {\it parabolic} if it contains $\frak g_0\oplus\frak g_+$ as a
proper subalgebra, that is
$$ 
\frak p = \bigoplus_{j\in \Bbb Z} \frak p_j, \quad\text{where $ \frak
p_j=\frak g_j$ for $j\ge 0$}, \hbox{ and } \frak p_j\neq 0 \hbox{ for  some
} j<0.
$$

\noindent We assume the following properties of $\frak g$:

(P1) $\frak g_0$ is commutative,

(P2) if $a\in \frak g_{-k}$ $ (k>0)$ and $[a, \frak g_1 ] =0$, then $a=0$.

\proclaim{Lemma 2.1} For any parabolic subalgebra $\frak p$ of $\frak g$,
$\frak p_{-k}\neq 0$, $k>0$,  implies  $\frak p_{-k+1}\neq 0$.
\endproclaim

\demo{Proof} If  $\frak p_{-k+1} = 0$, then $[\frak p_{-k} , \frak g_1]=0$,
i.e. for all $a\in \frak p_{-k}$, $[a,\frak g_1]=0$, and using (P2), we get
$a=0$.
\qed\enddemo

Given $a\in \frak g_{-1}$, $a\neq 0$, we define $\frak p^a = \oplus_{j\in
\Bbb Z} \frak p_j^a$, where $\frak p^a_j=\frak g_j$ for all $j\geq 0$, and
$$
\frak p^a_{-1}= \sum [...[[a,\frak g_0], \frak g_0],...]\ , \quad\quad
\frak p^a_{-k-1}=[\frak p^a_{-1}, \frak p^a_{-k}].
$$

\proclaim{Lemma 2.2} a) $\frak p^a$ is the minimal parabolic subalgebra
containing $a$.

\noindent b) $\frak g^a_0 := [\frak p^a , \frak p^a]\cap \frak g_0 = [ a,
\frak g_1]$.
\endproclaim

\demo{Proof} a) We  have to prove that $\frak p^a$ is a subalgebra. First,
$[\frak p^a_{-k}, \frak p^a_{-l}]\subseteq \frak p^a_{-l-k}$ $(k,l>0)$ is
proved by induction on $k$:
$$
\aligned [\frak p^a_{-k}, \frak p^a_{-l}] &= [[(\frak p^a_{-1})^{k-1} ,
\frak p^a_{-1} ] , (\frak p^a_{-1})^l]\\
&= [[(\frak p^a_{-1})^{k-1} ,  (\frak p^a_{-1})^l], \frak p^a_{-1}] +
[(\frak p^a_{-1})^{k-1} , [\frak p^a_{-1} , (\frak p^a_{-1})^l]]\\
&\subseteq [(\frak p^a_{-1})^{l+k-1} , \frak p^a_{-1} ] + [(\frak
p^a_{-1})^{k-1} , (\frak p^a_{-1})^{l+1}]\\
&\subseteq (\frak p^a_{-1})^{k+l}.
\endaligned
$$

And $[\frak p^a_{-k}, \frak g_{m}]\subseteq \frak p^a_{m-k}$  ($m<k$) also
follows by induction on $k$:
$$
\aligned [\frak p^a_{-k}, \frak g_{m}]&= [[(\frak p^a_{-1})^{k-1} , \frak
p^a_{-1}] , \frak g_m]\\
&= [[(\frak p^a_{-1})^{k-1} ,  \frak g_m],  \frak p^a_{-1}] + [ (\frak
p^a_{-1})^{k-1} , [\frak p^a_{-1} , \frak g_m]]\\
&\subseteq [ \frak p^a_{m-k+1}, \frak p^a_{-1} ] + [   (\frak
p^a_{-1})^{k-1} , \frak g_{m-1}]\\
&\subseteq \frak p^a_{m-k}.
\endaligned
$$
Finally, it is obviously the minimal one, proving (a).

b) For any $k>1$:
$$
\aligned [\frak p^a_{-k}, \frak g_{k}]&= [[(\frak p^a_{-1})^{k-1} ,  \frak
g_k],  \frak p^a_{-1}] + [ (\frak p^a_{-1})^{k-1} , [\frak p^a_{-1} , \frak
g_k]]\\
&\subseteq [  \frak g_1, \frak p^a_{-1} ] + [   (\frak p^a_{-1})^{k-1} ,
\frak g_{k-1}]
\endaligned
$$
Therefore, by induction, $\frak g^a_0 =  [  \frak g_1, \frak p^a_{-1} ] $. But
$$
\aligned
[  \frak g_1, \frak p^a_{-1} ]&= \hbox{ linear span } \{ [...[[a, c_1],
c_2],...], x] : c_i\in \frak g_0 , x\in \frak g_1 \}\quad \hbox{ (using
(P1))}\\
&=  \hbox{ linear span } \{ [a, [ c_1, ... [c_{k-1},[ c_k , x]...] : c_i\in
\frak g_0 , x\in \frak g_1 \}\\
&=[a, \frak g_1].
\endaligned
$$
proving the lemma. \qed
\enddemo

In the particular case of the central extension of the Lie algebra of
matrix differential operators on the circle (see [BKLY], Remark 2.2), we
observed the existence of some parabolic subalgebras $\frak p$ such that
$\frak p_{-j} = 0$ for $j>>0$. Having in mind that example, we give the
following

\definition{Definition 2.3} a) A parabolic subalgebra $\frak p$ is called
{\it non-degenerate} if $ \frak p_{-j} $ has finite codimension in $\frak
g_{-j}$, for all $j>0$.

\noindent b) An element $a\in \frak g_{-1}$ is called {\it non-degenerate}
if $\frak p^a$ is non-degenerate.

\enddefinition

Now, we begin our study of quasifinite representations over $\frak g$. A
$\frak g$-module $V$ is called $\Bbb Z${\it -graded} if $V=\oplus_{j\in
\Bbb Z} V_j$ and $\frak g_i V_j\subset V_{i+j}$. A $\Bbb Z\,  $-graded
$\frak g$-module $V$ is called {\it quasifinite} if $\dim V_j <\infty$ for
all $j$.

Given $\lambda\in \frak g_0^*$, a {\it highest weight module}  is a $\Bbb
Z$-graded $\frak g$-module $V(\frak g, \lambda )$ generated by a highest
weight vector $v_{\lambda }\in V(\frak g, \lambda )_0$ which satisfies
$$
h v_{\lambda } = \lambda (h) v_{\lambda } \qquad (h\in \frak g_0 ), \qquad
\frak g_+ v_{\lambda } =0.
$$
A non-zero vector $v\in V(\frak g, \lambda )$ is called {\it singular} if
$\frak g_+ v =0$. 

The  {\it Verma module} over $\frak g$ is defined as usual:
$$
M(\frak g, \lambda)= \Cal U(\frak g )\otimes_{\Cal U(\frak g_0 \oplus \frak
g_+)} \Bbb C_{\lambda}
$$
where $\Bbb C_{\lambda}$ is the 1-dimensional $(\frak g_0 \oplus \frak
g_+)$-module given by $h\mapsto \lambda(h)$ if $h\in \frak g_0$,  $\frak
g_+\mapsto 0$, and the action of $\frak g$ is induced by the left
multiplication in $\Cal U(\frak g )$. Here and further $\Cal U(\frak q)$
stands for the universal enveloping algebra of the Lie algebra $\frak q$.
Any  highest weight module $V(\frak g, \lambda )$ is a quotient module of
$M(\frak g, \lambda)$. The irreducible module $L(\frak g, \lambda)$ is the
quotient of $M(\frak g, \lambda)$ by the maximal proper graded submodule.
We shall write $M( \lambda)$ and $L( \lambda)$ in place of $M(\frak g,
\lambda)$ and $L(\frak g, \lambda)$ if no ambiguity may arise.

 Consider a parabolic subalgebra $\frak p=\oplus_{j\in \Bbb Z} \frak p_j$
of $\frak g$ and let $\lambda \in \frak g_0^*$ be such that $\lambda
|_{\frak g_0\cap [\frak p, \frak p]} =0$. Then the $(\frak g_0 \oplus \frak
g_+)$-module $\Bbb C_{\lambda}$ extends to a $\frak p$-module by letting
$\frak p_j$ act as 0 for $j<0$, and we may construct the highest weight
module 
$$
M(\frak g,\frak p, \lambda)= \Cal U(\frak g )\otimes_{\Cal U(\frak p)} \Bbb
C_{\lambda}
$$
called the  {\it generalized Verma module}.

\ 

We will also require the following condition on $\frak g$:

\vskip .3cm

(P3) If $\frak p$ is a non-degenerate parabolic subalgebra of $\frak g$,
then there exists a non-degenerate element $a$ such that $\frak
p^a\subseteq\frak p$.

\vskip .4cm

\remark{Remark 2.4} In all the examples considered in [KR1], [BKLY], [KWY]
and Section 3 of this work, property (P3) is satisfied.
\endremark

\proclaim {Theorem 2.5} The following conditions  on $\lambda \in \frak
g_0^* $ are equivalent:
\roster
\item $M(\lambda)$ contains   a singular vector $a . v_{\lambda}$ in
$M(\lambda)_{-1}$, where $a$ is non-degenerate; 
\item There exist a non-degenerate element $a\in \frak g_{-1}$, such that
$\lambda([\frak g_1 , a]) = 0$.
\item $L(\lambda)$ is quasifinite; 
\item There exist a non-degenerate element $a\in \frak g_{-1}$, such that
$L(\lambda)$ is the irreducible  quotient of the generalized Verma module
$M(\frak g , \frak p^a,\lambda)$. \endroster
\endproclaim

\demo{Proof} $(1)\Rightarrow (4):$ Denote by $a\, v_{\lambda}$ the singular
vector, where $a\in\frak g_{-1}$, then (4) holds for this particular $a$.
$(4)\Rightarrow (3)$ is immediate. Finally, $L(\lambda)$  quasifinite
implies $dim(\frak g_{-1}. v_{\lambda})< \infty$, then there exist an $a\in
\frak g_{-1}$ such that $a\, v_{\lambda}=0$ in $L(\lambda)$, so $0=\frak
g_{1}. (a v_{\lambda})= a (\frak g_{1}. v_{\lambda}) + [\frak g_{1}, a]
v_{\lambda}= \lambda([\frak g_1 , a])\  v_{\lambda}$, getting
$(3)\Rightarrow (2)\Rightarrow (1)$.
\qed
\enddemo

%
%
%
%

\head 3. The Lie algebra $W_{\infty}^{(n)}$ and its parabolic
subalgebras\endhead

\ 

We turn now to a certain family of  $\Bbb Z$-graded Lie algebras. 
Let   ${\Cal D} $ be the Lie algebra of  regular  differential operators on
the circle, i.e. the operators on $\Bbb C[t, t^{-1}]$ of the form 
$$
E=e_k(t)\partial_t^k + e_{k-1}(t)\partial_t^{k-1} + \dotsb +e_0(t), \hbox{
where } e_i(t)\in \Bbb C[t,t^{-1}],
$$
The elements 
$$
J^l_k= - t^{l+k} (\partial_t)^l \qquad (l\in \Bbb Z_+ , \  k\in \Bbb Z )
$$
form its  basis, where $\partial_t$ denotes $\frac d {dt}$. Another basis
of ${\Cal D} $ is 
$$
L^l_k = - t^k D^l \qquad (l\in \Bbb Z_+ , \  k\in \Bbb Z )
$$
where $D= t \partial_t$. It is easy to see that 
$$
J^l_k = - t^k [D]_l.\eqno (3.1)
$$
Here and further we use the notation 
$$
[x]_l = x (x-1) \dots (x-l+1).
$$

Fix a  linear map $\hbox{T}:\Bbb C [w] \to \Bbb C$. Then we have the
following 2-cocycle on ${\Cal D}$, where $f(w), g(w) \in \Bbb C [w]$ [KR1]:
$$
\Psi_{\hbox{T}} (z^rf(D),z^s g(D))=\cases  \hbox{T}(\sum\limits_{-r\le m\le
-1} f(w+m)g(w+m+r)) &\text{\rm{if} $r=-s\ge 0$},\\
0 &\text{\rm{if} $r+s\neq 0$.}\endcases \eqno (3.2)
$$

We let $\Psi= \Psi_{\hbox{T}}$ if $\hbox{T}: \Bbb C [w] \to \Bbb C$ is the
evaluation map at $w=0$. The central extension of ${\Cal D}$ by a
one-dimensional center $\Bbb C C$, corresponding to the 2-cocycle $\Psi$ is
 denoted  by $W_{1+\infty}$. The bracket in $W_{1+\infty}$ is given by 
$$
[t^rf(D) , t^sg(D)]= t^{r+s}(f(D+s)g(D)- f(D) g(D+r)) + \Psi (t^rf(D),t^s
g(D)) C. \eqno (3.3)
$$

Consider the following family of Lie subalgebras of ${\Cal D}$ ($n\in \Bbb
N$):
$$
\Cal D^{(n)} := \Cal D \partial_t^n.
$$

Denote by $W_{\infty}^{(n)}$ the central extension of $\Cal D^{(n)}$ by
$\Bbb C C$ corresponding to the restriction of the 2-cocycle $\Psi$.
Observe that $W_{\infty}^{(1)}$ is the well known $W_{\infty}$ subalgebra
of $W_{1+\infty}$.

Letting wt $t^kf(D) = k$, wt $C=0$ defines the {\it principal gradation} of
$W_{1+\infty}$ and of $W_{\infty}^{(n)}$:
$$
W_{\infty}^{(n)}= \bigoplus_{j\in \Bbb Z} (W_{\infty}^{(n)})_j,
$$
where using (3.1), we have explicitly:
$$
(W_{\infty}^{(n)})_j = \{ t^j f(D) [D]_n \  |  \  f(w)\in \Bbb C [w] \} +
\delta_{j0} \Bbb C.
$$

It is easy to check that the $\Bbb Z$-graded Lie algebras
$W_{\infty}^{(n)}$ satisfy the properties (P1-2).

\remark{Remark 3.4} The Lie algebra $W_{\infty}^{(n)}$ contains a $\Bbb
Z$-graded subalgebra isomorphic to the Virasoro algebra if and only if
$n=0$ or $n=1$. Indeed, from the commutator: 
$
[tf(D)[D]_n , g(D)[D]_n] = t\left(g(D)(D-n+1)-g(D+1)(D+1)\right) [D]_{n-1}
f(D) [D]_n
$, 
it is immediate that if the relation $[L_1, L_0]=L_1$ is satisfied (for
some elements $L_i\in (W_{\infty}^{(n)})_{i}$, $i=0,1$), then $n\leq 1$.
The existence of Virasoro subalgebras for $n\leq 1$ was observed in [KR1].
\endremark

\

Let  $\frak p$ be a parabolic subalgebra of $W_{\infty}^{(n)}$. Observe
that  for each $j\in \Bbb N$ we have: 
$$
\frak p_{-j} =\{t^{-j} f(D) \   | \  f(w)\in I_{-j} \}
$$ 
where $I_{-j}$ is a subspace of $[w]_n \Bbb C[w] $. Since
$$
[f(D)[D]_n, t^{-k}p(D)]=t^{-k}(f(D-k)[D-k]_n  - f(D)[D]_n) p(D)
$$
we see that $I_{-k}$ satisfies 
$$
A_{n,k}\  . \ I_{-k} \subseteq I_{-k}
$$
where $A_{n,k}=\{ f(w-k)[w-k]_n  - f(w)[w]_n  \  | \   f(w)\in \Bbb C[w] \}$. 

\proclaim{Lemma 3.5} a) $I_{-k}$ is an ideal for all $k\in \Bbb N$  if
$n\leq 2$ (there are examples of parabolic subalgebras where $I_{-1}$ is
not an ideal for any $n>2$).

\noindent b) If $I_{-k}\neq 0$, then it has finite codimension in $\Bbb C
[x]$.
\endproclaim

\demo{Proof} Observe that $A_{0,k}=A_{1,k}= \Bbb C [w]$ for all $k\geq 1$,
and $A_{2,k}$ is a subspace which contains a polynomial of degree $l$ for
all $l\geq 1$, proving the first part. Now, observe that $A_{n,1}= \Bbb
C[w] [w-1]_{n-1}$, and take 
$$
I_{-1}= \Bbb C \ [w]_n \oplus \Bbb C\ w[w]_n\oplus A_{n,1} [w]_n \oplus
A_{n,1} w[w]_n,
$$
$$
I_{-k}=\Bbb C[w] [w-k+1]_n\dots[w-1]_n [w]_n, \qquad k>1.
$$ 
Then, after some computation, it is possible  to see that these subspaces
define a parabolic subalgebra $\frak p= \oplus_{k\in\Bbb Z}\frak p_k $,
where $\frak p_{-k}=\{t^{-j} f(D) \   | \  f(w)\in I_{-j} \}$ for $k\geq
1$. Observe that for $n>2$, $I_{-1}$ is no longer an ideal. 

Finally, since $A_{n,k}$ contains a polynomial of degree $l$ for all $l\geq
m$, for some $m\in \Bbb N$, part b) follows. \qed\enddemo

\vskip .2cm

\remark{Remark 3.6} Due to Lemma 3.5, we no longer have the situation of
parabolic subalgebras described in terms of ideals, as in the Lie algebras
considered in [KR1], [BKLY] and [KWY]. But for the algebra  $W_{\infty}$
the parabolic subalgebras are as in these references. 
\endremark

\

We shall need the following proposition to  study modules over
$W_{\infty}^{(n)} $ induced from its parabolic subalgebras.

\proclaim {Proposition 3.7}  a) Any non-zero element $d\in
(W_{\infty}^{(n)})_{-1}$ is non-degenerate.

\noindent b) Any parabolic subalgebra of $W_{\infty}^{(n)}$ is non-degenerate.

\noindent c) Let $d=t^{-1}b(D)=t^{-1}a(D)[D]_n \in
(W_{\infty}^{(n)})_{-1}$, then
$$
\aligned
(W_{\infty}^{(n)})_0^d :&= [(W_{\infty}^{(n)})_1 , d] \\
&= \hbox{span}\{[D-1]_ng(D)b(D)  - [D]_n g(D+1) b(D+1) \ |\  g \in \Bbb
C[w]\}.
\endaligned
$$
\endproclaim

\demo{Proof}  Let $0\neq d\in (W_{\infty}^{(n)})_{-1}$, then $\frak
p^d_{-j}\neq 0$ for all $j\geq 1$. So, by Lemma 3.5 b), part a) follows .
Let $\frak p$ be any parabolic subalgebra of $W_{\infty}^{(n)}$,  using
Lemma 2.1 we get $\frak p_{-1}\neq 0$. Then, using a) and $\frak p^d
\subseteq \frak p$ (for any non-zero $d\in \frak p_{-1}$), we obtain b).
Finally, part c) follows by Lemma 2.2 b) and the commutator: 
$
[t g(D) [D]_n , t^{-1}b(D)] = [D-1]_n g(D-1) b(D) - [D]_n g(D) b(D+1).
$
\qed\enddemo

%
%
%
%

\head 4. Quasifinite Highest Weight Modules over $W_{\infty}^{(n)} $ \endhead

\

By Proposition 3.7,  $W_{\infty}^{(n)} $ also satisfies property (P3),
hence we can apply Theorem 2.5.

Let $L(\lambda)$ be a quasifinite highest weight module over
$W_{\infty}^{(n)}$. By Theorem 2.5, there exists some monic polynomial
$b(w)=a(w)  [w]_n$ such that
$$
(t^{-1}b(D)) v_{\lambda}=0 . 
$$
We shall call such  monic polynomial of minimal degree, uniquely determined
by the highest weight $\lambda$, the {\it characteristic} polynomial of
$L(\lambda)$. 

\vskip .3cm

A functional $\lambda \in (W_{\infty}^{(n)})_0^* $ is described by its {\it
labels} $\Delta_{l}=-\lambda(D^l \, [D]_n)$, where $l\in \Bbb Z_+$, and the
{\it central charge} $c=\lambda(C)$. We shall consider  the generating series
$$
\Delta_{\lambda }(x)=\sum_{l=0}^{\infty} \frac{x^l}{l!} \Delta_{l}. \eqno
(4.1)
$$

Recall that a {\it quasipolynomial} is a linear combination of functions of
the form $p(x)e^{\alpha x}$, where $p(x)$ is a polynomial and $\alpha\in
\Bbb C$. Recall the well-known  characterization: a formal power series is
a quasipolynomial if and only if it satisfies a non-trivial linear
differential equation with constant coefficients. We have the following
characterization of quasifinite highest weight modules.

\proclaim {Theorem 4.2} A $W_{\infty}^{(n)} $-module $L(\lambda)$ is
quasifinite if and only if 
$$
n \Delta_{\lambda}(x) + \frac d {dx} ( (e^x - 1) \Delta_{\lambda}(x))
$$
is a quasipolynomial. This condition is equivalent to the existence of  a
quasipolynomial $\phi_{\lambda} (x)$ with  $\phi_{\lambda} (0)=0$,  such
that  
$$
\Delta_{\lambda}(x)= \left[\frac d {dx}\right]_n \left(\frac
{\phi_{\lambda}(x)}{e^x - 1}\right).\eqno (4.3) 
$$
\endproclaim

\demo{Proof} From Proposition 3.7 a) and c),  and Theorem 2.5(2), we have
that $L(\lambda)$ is quasifinite if and only if  there exist a polynomial 
$
b(w)=[w]_n a(w)
$ 
such that 
$$
\lambda\left([D-1]_ng(D)b(D)  - [D]_n g(D+1) b(D+1)\right) = 0 \quad \hbox{
for any polynomial } g,
$$
or equivalently
$$
0=\lambda\left(   b(D) [D-1]_n e^{x D} -  b(D+1) [D]_n  e^{x(D+1)} \right).
\eqno(4.4)
$$
Now, using that $\Delta_{\lambda}(x)= -\lambda \left(\left[\frac d
{dx}\right]_n   e^{x D} \right)$, and the identities 
$$
f(D) e^{x D} = f\left(\frac d {dx}\right) \left(e^{x D}\right),\quad \quad
\quad [D]_n  e^{x(D+1)}= e^x [D]_n  e^{xD} = e^x  \left[\frac d
{dx}\right]_n e^{x D},
$$
condition (4.4) can be rewritten as follows:
$$
\aligned
0&=\lambda\left(   b(D) [D-1]_n e^{x D} -  b(D+1) [D]_n  e^{x(D+1)} \right)\\
&=\lambda\left(  b\left(\frac d {dx}\right) \left[\frac d {dx} - 1
\right]_n (e^{x D}) -  b\left(\frac d {dx}\right) \left( [D]_n
e^{x(D+1)}\right)\right)\\
&= \lambda\left(   a\left(\frac d {dx}\right) \left[\frac d {dx} - 1
\right]_n \left[\frac d {dx}\right]_n (e^{x D}) -  b\left(\frac d
{dx}\right) \left(e^x  \left[\frac d {dx}\right]_n e^{x D}\right)\right)\\
&= - a\left(\frac d {dx}\right) \left[\frac d {dx} - 1 \right]_n
\Delta_{\lambda}(x)  +  b\left(\frac d {dx}\right) (e^x \Delta_{\lambda}(x))\\
&= - a\left(\frac d {dx}\right) \left[\frac d {dx} - 1 \right]_{n-1}
\left(\left(\frac d {dx} - n\right) \Delta_{\lambda}(x)  +  \left(\frac d
{dx}\right) (e^x \Delta_{\lambda}(x))\right).
\endaligned
$$
Thus, $L(\lambda)$ is quasifinite if and only if there exists a polynomial
$a(w)$ such that 
$$
a\left(\frac d {dx}\right) \left[\frac d {dx} - 1 \right]_{n-1} \left(n
\Delta_{\lambda}(x) + \frac d {dx} ( (e^x - 1) \Delta_{\lambda}(x))\right)
= 0 
$$
Therefore, $L(\lambda)$ is quasifinite if and only if $F_{\lambda}(x):= n
\Delta_{\lambda}(x) + \frac d {dx} ( (e^x - 1) \Delta_{\lambda}(x))$ is a
quasipolynomial, proving the first equivalence of the theorem.

Now, suppose that $\Delta_{\lambda}(x)= \left[\frac d {dx}\right]_n
\left(\frac {\phi_{\lambda}(x)}{e^x - 1}\right)$ with $\phi_{\lambda}
(0)=0$, then using that  
$$
\left[\frac d {dx} - 1 \right]_n ( e^x g(x)) \, = \, e^x \left[\frac d
{dx}\right]_n (g(x))
$$
we have 
$$
F_{\lambda}(x)=\left[\frac d {dx}\right]_{n+1} \phi_{\lambda}(x).\eqno (4.5)
$$
Thus, if $\phi_{\lambda}(x)$ is a quasipolynomial, then $F_{\lambda}(x)$ is
a quasipolynomial.

Conversely, suppose that $F_{\lambda}(x)$ is a quasipolynomial. Let
$\phi(x)$ be {\bf any} solution of (4.5) satisfying the condition $\phi
(0)=0$. Observe that all solutions are of the form $\phi(x) + \sum_{i=1}^n
c_i (e^x - 1)^i$.  Then by (4.5), $\phi(x)$ is a quasipolynomial, and we
have $\Delta_{\lambda}(x)= \left[\frac d {dx}\right]_n \left(\frac
{\phi(x)}{e^x - 1}\right)$, where the expression of $ \Delta_{\lambda}(x)$
is independent of the choice of $\phi(x)$, completing the proof.   
\qed\enddemo

\remark{Remark 4.6} Given a quasipolynomial $F(x)$, it is easy to see that
$\Delta_{\lambda}(x)$ is uniquely determined by the equation $F (x) = n
\Delta_{\lambda}(x) + \frac d {dx} ( (e^x - 1) \Delta_{\lambda}(x))$.
\endremark

\definition{Definition 4.7} The quasipolynomial $\phi_{\lambda}(x)+c$,
where $\phi_{\lambda}(x)$ is from  (4.3) and $c$ is the central charge,
can be (uniquely) written  in the form 
$$
\phi_{\lambda}(x)+c=\sum_{r\in I} p_{r}(x) e^{rx} ,\tag 4.8
$$
where all $r$ are distinct numbers. 
The numbers $r$ appearing in (4.8) are called {\it exponents} of the
$W_{\infty}^{(n)} $-module $L(\lambda)$, and  the polynomial $p_{r}(x)$ is
called the {\it multiplicity} of $r$, denoted by mult($r$). Note that, by
definition:
$$
c=\sum_r p_r (0).
$$
\enddefinition

\proclaim{Corollary 4.9} Let $L(\Lambda)$ be a quasifinite irreducible
highest weight module over $W_{\infty}^{(n)} $, let $b(w)=[w]_n a(w)$ be
its  characteristic polynomial, and let $F_{\lambda}(x)=$\linebreak $ n
\Delta_{\lambda}(x)  +   \frac d {dx} ( (e^x - 1) \Delta_{\lambda}(x))$.
Then $\ a\left(\frac d {dx} \right)\left[\frac d {dx} -1 \right]_{n-1} F(x)
= 0$ is the    minimal order homogeneous linear differential equation with
constant coefficients of the form $f\left(\frac d {dx} \right)\left[\frac d
{dx} -1 \right]_{n-1}$, satisfied by $F(x)$. Moreover, the exponents
appearing in (4.8) are all   roots of the polynomial $[w-1]_{n-1}a(w)$.
\endproclaim

\ 

Now, we will consider the restriction of quasifinite highest weight modules
over $W_{1+\infty}$ to $W_{\infty}^{(n)}$. We will need some notation.

A functional $\lambda\in (W_{1+\infty})_0^*$ is characterized by its labels
$\Gamma_m=-\lambda(D^m)$, where $m\in\Bbb Z_+$, and the central charge
$c=\lambda(C)$, cf. (4.1). Introduce the new generating series:
$$
\Gamma_{\lambda}(x)=\sum_{m=0}^{\infty} \frac{x^m}{m!} \Gamma_m.
$$
Observe that $\left[\frac d {dx}\right]_n \Gamma_{\lambda}(x) =
\Delta_{\lambda}(x)$, where $\Delta_{\lambda}(x)$ is defined by (4.1).
Recall that a $W_{1+\infty}$-module $L(W_{1+\infty}, \lambda)$ is
quasifinite if and only if 
$
\Gamma_{\lambda}(x)=\frac {\phi (x)}{e^x-1}
$,  
where $\phi(x)$ is  a quasipolynomial such that $\phi(0)=0$ [KR1]. We have
the following partial restriction result:

\proclaim{Proposition 4.10} Any quasifinite  $W_{\infty}^{(n)}$-module
$L(\lambda)$ can be obtained as a quotient of the
$W_{\infty}^{(n)}$-submodule generated by the highest weight vector  of a
quasifinite $W_{1+\infty}$-module $L(W_{1+\infty}, \tilde{\lambda})$, for
some quasifinite  $\tilde{\lambda}\in (W_{1+\infty})_0^*$ such that
$\tilde{\lambda}_{|_{(W_{\infty}^{(n)})_0}} = \lambda$.

\endproclaim

\demo{Proof} Given $\Delta_{\lambda}(x)= \left[\frac d {dx}\right]_n
\left(\frac {\phi (x)}{e^x-1}\right)$, consider $\tilde{\lambda}\in
(W_{1+\infty})_0^*$ determined by  \linebreak $\Gamma_{\tilde{\lambda}}(x)
=  \frac {\phi (x)}{e^x-1}$, and  the proposition follows. \qed

\enddemo

\

Let ${\Cal O}$ be the algebra of all holomorphic functions on $\Bbb C$ with
the topology of uniform convergence on compact sets. We consider the vector
space ${\Cal D^{\Cal O}}$ spanned by the differential operators (of
infinite order) of the form $t^kf(D)$, where $f\in {\Cal O}$. The bracket
in ${\Cal D}$ extends to ${\Cal D^{\Cal O}}$. Then the cocycle $\Psi$
extends to a 2-cocycle on $ {\Cal D^{\Cal O}}$ by formula (3.2). Let
$W_{1+\infty}^{\Cal O}= {\Cal D^{\Cal O}} +\Bbb C C$ be the corresponding
central extension with the principal gradation as in $W_{1+\infty}$.

Consider the Lie subalgebras of $\Cal D^{\Cal O}$:
$$
\Cal D_{\Cal O}^{(n)}:= \Cal D^{\Cal O} \ \partial_t^n. \eqno (4.11)
$$
We shall denote by $W_{\infty}^{(n) \Cal O}$ the central extension of $\Cal
D_{\Cal O}^{(n)}$ by $\Bbb C C$ corresponding to the restriction of the
cocycle $\Psi$. And we shall use de notation $W_{\infty}^{\Cal O} =
W_{\infty}^{(1) \Cal O}$. Observe that $W_{\infty}^{(n) \Cal O}$ inherit a
$\Bbb Z\, $-gradation from $W_{1+\infty}^{\Cal O}$.

In the following section, we shall need the following proposition.

\proclaim{Proposition 4.12} Let $V$ be a quasifinite  $W_{\infty}^{(n)}
$-module. Then the action of $W_{\infty}^{(n)} $ on $V$ naturally extends
to the action of $(W_{\infty}^{(n) \Cal O})_k$ on $V$ for any $k\neq
0$.\endproclaim

\demo{Proof} The proof is analogous to that of  Proposition 4.3 in
[KR1].\qed\enddemo

%
%
%
%

\head 5. Embedding of $W_{\infty}$ into $\widehat{gl}_{\infty}$.
\endhead

\

In the following, {\bf we will suppose that $n=1$, i.e. we consider the
algebra $W_{\infty}$}.

Let $gl_{\infty}$ be the $\Bbb Z$-graded Lie  algebra of all matrices
$(a_{ij})_{i,j\in\Bbb Z}$ with finitely many nonzero diagonals (deg
$E_{ij}=j-i$). Consider the central extension
$\widehat{gl}_{\infty}=gl_{\infty} + \Bbb C C$ defined by the cocycle:
$$
\Phi(A,B)=\hbox{tr } ([J,A]B),\qquad J=\sum_{i\leq 0} E_{ii}.
$$

Given $s\in \Bbb C$, we will consider the natural action of the Lie algebra
$\Cal D^{(1)}$ (resp. $\Cal D_{\Cal O}^{(1)}$) on $t^s \Bbb C[t,t^{-1}]$.
Taking the basis $v_j=t^{-j+s}$ $(j\in\Bbb Z)$ of this space, we obtain a
homomorphism of  Lie algebras $\varphi_s : \Cal D^{(1)}\to gl_{\infty}$
   (resp.  $\varphi_s : \Cal D_{\Cal O}^{(1)}\to gl_{\infty}$):
$$
\varphi_s(t^kf(D) D)=\sum_{j\in\Bbb Z} f(-j+s)(-j+s) E_{j-k,j}.
$$
This homomorphism preserves gradation and it lifts to a homomorphism
$\widehat{\varphi}_s$ of the corresponding central extensions as follows
[KR1]:
$$
\widehat{\varphi}_s(De^{xD})= \varphi_s(De^{xD})- \left(\frac{e^{sx}
-1}{e^x -1}\right)'  C,\qquad \widehat{\varphi}_s(C)=C.
$$

Let $s\in \Bbb Z $ and denote by  $\widehat{gl}_{\infty,s}$  the Lie
subalgebra of $\widehat{gl}_{\infty}$ generated by $C$ and $\{E_{ij} |
i\neq -s \hbox{ and }  j\neq s \}$. Observe that $\widehat{gl}_{\infty,s}$
is naturally isomorphic to $\widehat{gl}_{\infty}$. Let
$p_s:\widehat{gl}_{\infty}\to \widehat{gl}_{\infty,s}$ be the projection
map. If $s\in \Bbb Z$, we  redefine $\widehat{\varphi}_s$ by the
homomorphism $p\circ \widehat{\varphi}_s:W_{\infty}\to
\widehat{gl}_{\infty,s}$.

Given $\vec s=(s_1,\dots , s_m)\in \Bbb C^m$, we have a homomorphism of Lie
algebras over $\Bbb C$:
$$
\widehat{\varphi}_{\vec s}=\oplus_{i=1}^m \widehat{\varphi}_{s_i} :
W_{\infty} \to \frak g_{\vec s}=\oplus_{i=1}^m \frak g_{s_i}\eqno (5.1)
$$
where $\frak g_{s_i}=\widehat{gl}_{\infty} $ if $s_i\notin\Bbb Z$, and
$\frak g_{s_i}=\widehat{gl}_{\infty,s_i} $ if $s_i\in \Bbb Z$. The proof of
the following proposition is similar to that of Proposition 3.2 in [KR1].

\proclaim{Proposition 5.2} The homomorphism $\widehat{\varphi}_{\vec s}$
extends to a  homomorphism of Lie algebras over $\Bbb C$, which is also
denoted by $\widehat{\varphi}_{\vec s}$:
$$
\widehat{\varphi}_{\vec s}:W_{\infty}^{\Cal O} \to \frak g_{\vec s}.
$$
The homomorphism $\widehat{\varphi}_{\vec s}$ is surjective provided that
$s_i-s_j\notin \Bbb Z$ for $i\neq j$.
\endproclaim

\remark{Remark 5.3} For $s\in \Bbb Z$ the image of $W_{\infty}^{\Cal O}$
under the homomorphism $\widehat{\varphi}_s$ is \linebreak $\nu^k
(\widehat{gl}_{\infty,s-k})$ for any $k\in\Bbb Z$, where $\nu$ is the
automorphism defined by 
$$
\nu(E_{ii})=E_{i+1,i+1}.\eqno (5.4)
$$
Hence we may (and will) assume that $0\leq \hbox{Re }s <1$ throughout the
paper.
\endremark

%
%
%
%

\head 6. Unitary Quasifinite Highest Weight Modules over $W_{\infty}$
\endhead

\  

The algebra $\Cal D^{(1)}$ acts on the space $V=\Bbb C[t,t^{-1}] / \Bbb C$.
One has a non-degenerate Hermitian form on $V$:
$$
B(f,g)= \hbox{Res}_t \  \overline f dg,
$$
where $(\overline{\sum a_i t^i})=\sum \overline{a_i} t^{-i}$, $a_i\in \Bbb
C$, (cf. [B]).

Consider the additive map $\omega :  \Cal D^{(1)} \to \Cal D^{(1)}$,
defined by:
$$
\omega(t^kf(D) D)=t^{-k}\overline{f}(D-k) D
$$
where for $f(D)=\sum_if_iD^i$, we let $\overline f(D)=\sum_i\overline
{f_i}D^i(f_i\in\Bbb C)$.

\proclaim{Proposition 6.1} The map $\omega$ is an anti-involution of the
Lie algebra $\Cal D^{(1)}$, i.e. $\omega$ is an additive map such that
$$
\omega^2 =id, \quad \omega(\lambda a)=\bar{\lambda} \omega(a), \ \ \hbox{
and }\  \omega([a,b])=[\omega(b),\omega(a)], \hbox{ for } \lambda\in \Bbb
C, a,b\in \Cal D^{(1)}.
$$
Furthermore, the operators $\omega(a)$ and $a$ are adjoint operators on $V$
with respect to $B$, and  $\omega ((\Cal D^{(1)})_j)=(\Cal D^{(1)})_{-j}$.
\endproclaim

\demo{Proof} The properties $\omega^2 =id, \  \omega(\lambda
a)=\bar{\lambda} \omega(a)$ are obvious. Now,
$$
\aligned
\omega([t^k &f(D)D, t^l g(D) D])= \\
&= \omega( t^{k+l} (f(D+l)(D+l)g(D)-g(D+k)(D+k)f(D))D)\\
&=  t^{-(k+l)} (\bar f(D-k)(D-k)\bar g(D-k-l)-\bar g(D-l)(D-l)\bar f(D-k-l))D.
\endaligned
$$
On the other hand,
$$
\aligned
[\omega(t^l &g(D)D), \omega(t^k f(D) D)]= [ t^{-l} \bar g(D-l) D,
t^{-k}\overline{f}(D-k) D]\\
&= t^{-(k+l)} (\bar g(D-k-l)(D-k)\bar f(D-k)-\bar f(D-k-l)(D-l)\bar g(D-l)) D.
\endaligned
$$
Hence $\omega$ is an anti-involution. Now we compute
$$
B(t^k f(D)D \ t^l , t^n)= B(l\, f(l)\,  t^{k+l}, t^n)=  n\, l\, \overline f
(l) \delta_{k+l,n}
$$
And we also have
$$
\aligned
B(t^l,\omega(t^kf(D)D)\  t^n)&= B(t^l, t^{-k}\overline f (D-k)D\  t^n)\\
&= B(t^l,n \, \overline f (n-k) t^{n-k})= n\, l\,  \overline f (n-k)
\delta_{l,n-k},
\endaligned
$$
proving the proposition.
\qed\enddemo

This anti-involution $\omega$ extends to the whole algebra $ {\Cal D_{\Cal
O}^{(1)}}$, defined in (4.11). Observe  that 
$$
\Psi(\omega(A),\omega(B))=\omega(\Psi(B,A)), \qquad A,B\in  {\Cal D_{\Cal
O}^{(1)}}.
$$
Therefore, the anti-involution $\omega$ of the Lie algebras $ \Cal D^{(1)}$
and $ {\Cal D_{\Cal O}^{(1)}}$ lifts to an anti-involution of their central
extensions $ W_{\infty}$ and  $W_{\infty}^{\Cal O}$, such that
$\omega(C)=C$, which we again denote by $\omega$.

In this section we shall classify and construct all  {\it unitary}
(irreducible) quasifinite highest weight modules over $W_{\infty}$ with
respect to the anti-involution  $\omega$. In order to do it, we shall need
the following lemma.

\proclaim{Lemma 6.2} Let $V$ be a unitary quasifinite highest weight module
over $W_{\infty}$ and let $b(w)=w a(w)$ be its first   characteristic
polynomial. Then $a(w)$ has only simple real roots.
\endproclaim

\demo{Proof} Let $v_{\lambda}$ be a highest weight vector of $V$. Then the
first graded subspace $V_{-1}$ has a basis
$$
\{ (t^{-1} D^{j+1}) v_{\lambda} \  | \  0\leq j < \hbox{deg } a \}.
$$
Consider the action of 
$$
S= -\frac 1 2 \left( D^2 + \left(\frac{1-\Delta_1}{1+\Delta_0} \right)
D\right)
$$
on $V_{-1}$. It is straightforward to check that 
$$
S^j ((t^{-1} D) v_{\lambda}) = (t^{-1} D^{j+1}) v_{\lambda}\qquad \hbox{for
all } j\geq 0.
$$
It follows that $a(S) ((t^{-1} D) v_{\lambda}) = 0$, and that $\{ S^j
((t^{-1} D) v_{\lambda}) \  | \  0\leq j <\hbox{deg }a   \}$  is a basis of
$V_{-1}$. We conclude from the above that $a(w)$ is the characteristic
polynomial of the operator $S$ on  $V_{-1}$. Since the operator $S$ is
self-adjoint, all the roots of $a(w)$ are real.

Now, suppose that $a(w)= (w-r)^m c(w)$ for some polynomial $c(w)$ and $r\in
\Bbb R$. Then $v= (S-r)^{m-1} c(S)((t^{-1} D) v_{\lambda})$ is a non-zero
vector in $V_{-1}$, but 
$$
(v ,v) = (c(S) ((t^{-1} D) v_{\lambda}), (S-r)^{2m-2} c(S)((t^{-1} D)
v_{\lambda})) = 0 \qquad \hbox{if } m\geq 2.
$$
Hence the unitarity condition forces $m=1$.
\qed\enddemo

\

Take $0<s<1$. Then  under the homomorphism $\widehat{\varphi}_s:
W_{\infty}^{\Cal O}  \to \widehat{gl}_{\infty}$, the anti-involution
$\omega$ induces the following   anti-involution on $\widehat{gl}_{\infty}$:

$$
\omega '(E_{ij}) = \left(\frac {s-i}{s-j}\right) E_{ji}, \qquad \omega '(C)=C.
$$
Indeed: 
$$
\aligned
\varphi_s(\omega(t^k f(D) D)) &= \varphi_s(t^{-k}\bar f (D-k) D)\\
&=\sum_{j\in \Bbb Z} \bar f (s-j-k) (s-j) E_{j+k,j}\\
&=\sum_{j\in \Bbb Z} \bar f (s-j) (s-j+k) E_{j,j-k}\\
&=\omega ' (\sum_{j\in \Bbb Z}  f (s-j) (s-j) E_{j-k,j}) \\
&=\omega ' (\varphi_s(t^k f(D) D)).
\endaligned
$$

Let $e_i=E_{i,i+1}$ and $f_i=E_{i+1,i}$, then $\omega' (e_i) = \lambda_i
f_i$ with $\lambda_i=\frac {s-i}{s-i-1}$. Observe that 
$$
\lambda_i<0 \quad \hbox{ if and only if }\quad i<s<i+1. \eqno (6.3)
$$
If we consider the linear automorphism $T$ defined by $T(e_i)=\mu_i
e_i=e'_i, \ T(f_i)=\mu_i^{-1} f_i =f'_i$ for some $\mu_i\in \Bbb C$, then
$\omega'(e'_i)=\omega(\mu_i e_i)=\bar\mu_i \lambda_i f_i= |\mu_i |^2
\lambda_i f'_i$. Hence by  (6.3),  $\omega'$ is equivalent to the
anti-involution $\tilde\omega$ defined by:
$$
\tilde{\omega}(e_i)=\cases f_i & \hbox{if } i\neq 0,\\
-f_i & \hbox{if } i= 0.
\endcases
$$
After a shift by  the automorphism  $\nu$  defined in (5.4), we may assume
$0<s<1$, then under the homomorphism $\widehat{\varphi}_s:W_{\infty}^{\Cal
O} \to \widehat{gl}_{\infty}$, the anti-involution $\omega$ induces an
anti-involution on $\widehat{gl}_{\infty}$ that is equivalent to the
following:
$$
E_{ij}^{\dag} = E_{ji}\  \hbox{ if }\  i,j>0 \ \hbox{ or }\  i,j\leq 0,
\quad E_{ij}^{\dag} = - E_{ji} \  \hbox{ otherwise}, \ \hbox{ and } \
C^{\dag}=C.
$$

As usual, for any $\lambda\in (\widehat{gl}_{\infty})_0^*$ we have the
associated irreducible highest weight $\widehat{gl}_{\infty}$-module
$L(\widehat{gl}_{\infty} , \lambda)$. An element $\lambda\in
(\widehat{gl}_{\infty})_0^*$ is determined by its {\it labels}
$\lambda_i=\lambda(E_{ii}), \ i\in \Bbb Z$, and {\it central charge}
$c=\lambda(C)$. Let $n_i=\lambda_i - \lambda_{i+1} +\delta_{i,0} c$
$(i\in\Bbb Z)$. The following classification is taken from [KR2] and [O]:

\proclaim{Proposition 6.4}  A non-trivial highest weight
$\widehat{gl}_{\infty}$-module with highest weight $\lambda$ and central
charge $c$  is unitary with respect to $^{\dag}$ if and only if  the
following properties hold:
$$
n_i\in \Bbb Z_+ \hbox{ if } i\neq 0 \hbox{ and } c=\sum_i n_i <0,\eqno (6.5a)
$$
$$
\hbox{if } n_i\neq 0 \hbox{ and } n_j\neq 0 , \hbox{ then } |i-j|\leq
-c.\eqno (6.5b)
$$
\endproclaim

In the case of $s=0$, the homomorphism $\widehat{\varphi}_0:
W_{\infty}^{\Cal O}  \to \widehat{gl}_{\infty,0}\simeq
\widehat{gl}_{\infty}$ induces the standard anti-involution on
$\widehat{gl}_{\infty}$: $(A)^*= \ ^t\bar A$. The following is a very well
known result:

\proclaim{Proposition 6.6}  A highest weight $\widehat{gl}_{\infty}$-module
with highest weight $\lambda$ and central charge $c$  is unitary with
respect to $^*$ if and only if  
$n_i\in \Bbb Z_+ $  and $c=\sum_i n_i $.
\endproclaim

Let $\lambda_i\in (\frak g_{s_i})_0^*$ such that $L(\frak g_{s_i},
\lambda_i)$ is a quasifinite $\frak g_{s_i}$-module. Then the tensor product
$$
 L(\frak g_{\vec s}, \vec \lambda):= \otimes_i L(\frak g_{s_i}, \lambda_i)
$$
is an irreducible $\frak g_{\vec s}\,  $-module.

\proclaim{Theorem 6.7} Let $V$ be a quasifinite  $\frak g_{\vec s}$-module,
viewed as a $W_{\infty}$-module via the homomorphism
$\widehat{\varphi}_{\vec s}$, where  $s_i - s_j\notin \Bbb Z$ if $i\neq j$,
and $0\leq \hbox{Re }s_i < 1$. Then any  $W_{\infty}$-submodule of $V$ is
also a $\frak g_{\vec s}$-submodule. In particular, the
$W_{\infty}$-modules $ L(\frak g_{\vec s}, \vec \lambda)$ are irreducible,
and in this way we  obtain all quasifinite $W_{\infty}$-modules
$L(\lambda)$ with $\phi_{\lambda}(x)=\sum_i n_i e^{r_i x}$, $n_i\in \Bbb
C$, $r_i\in \Bbb C$.
\endproclaim 

\demo{Proof} Consider any $W_{\infty}$-submodule $W$ of $V$. By Proposition
4.12, the action of $W_{\infty}$ can be extended to $(W_{\infty}^{\Cal
O})_k$ ($k\neq 0$). Using Proposition 5.2, we see that the subspace $W$ is
preserved by $\frak g_{\vec s} $. Therefore, the  $W_{\infty}$-modules $
L(\frak g_{s},  \lambda)$ are quasifinite and irreducible. Then it is easy
to calculate the generating series of the highest weight (see Section 4.6
in [KR1]): in the case  $s\in\Bbb R\backslash \Bbb Z$, we have
$$
\Delta_{s,\lambda}(x)= - \lambda(\hat{\varphi}_s (D e^{xD}))=\frac d {dx}
\left(\frac {\sum_{k\in\Bbb Z} e^{(s-k)x} n_k - c}{e^x -1}\right),\eqno (6.8)
$$
in the case of $s=0$, we have
$$
\Delta_{0,\lambda}(x)=\frac d {dx} \left(\frac {\sum_{j>0} e^{-j x} n_j +
\sum_{j<0} e^{(-j+1) x} n_j  + e^{2x} \lambda_0 - e^{-x} \lambda_1}{e^x
-1}\right),\eqno (6.9)
$$
and the last part of the theorem follows from equations (6.8-9). \qed\enddemo

\

In fact, as in Theorem 4.6 in [KR1], it is possible to construct all
irreducible quasifinite $W_{\infty}$-module in terms of representations of
$\widehat{gl}(\infty, R_m)$ (or a subalgebra of it), where
$\widehat{gl}(\infty, R_m)$ is the central extension of the Lie algebra of
infinite matrices with finitely many non-zero diagonals and coefficients in
the algebra of truncated polynomials $R_m:= \Bbb C[u] /(u^{m+1})$. More
precisely, we consider the homomorphism $\varphi_s^{[m]} : {\Cal D}^{(1)}
\to gl(\infty, R_m)$ given by 
$$
\varphi_s^{[m]} (t^k f(D) D)=\sum_{i=0}^m \sum_{j\in \Bbb Z} \frac{(f^{(i)}
(s-j) (s-j) + i f^{(i-1)}(s-j))}{i!} u^i E_{j-k,j}.
$$
In the case $s\in \Bbb R\backslash \Bbb Z$ we take $\frak g^{[m]}_s =
\widehat{gl}(\infty, R_m)$, and for $s\in \Bbb Z$ we have to remove the
generators $E_{r,s}$ and $u^m E_{-s,r}$ for all $r\in \Bbb Z$. All
quasifinite irreducible $L(\lambda)$ can be obtained using representations
of the Lie algebra $\frak g^{[\vec m]}_{\vec s}=\oplus_i \frak
g^{[m_i]}_{s_i}$ via the homomorphism $\varphi_{\vec s}^{[\vec m]}=\oplus_i
\varphi_{s_i}^{[m_i]}$, and as in [KR1] the coefficients in $\vec m$ are
given by the degree of the (polynomial) multiplicities of $\phi_{\lambda}(x)$.

\proclaim{Lemma 6.10} Only those highest weight representations of
$\widehat{gl}(\infty, R_m)$ that factor through $\widehat{gl}(\infty, \Bbb
C)$ are unitary.
\endproclaim

\demo{Proof} Indeed, let $v$ be a highest weight vector. Fix $i\in \Bbb Z$
and let $e=E_{i,i+1},\   f=E_{i+1,i}, \  h=E_{ii}$. Now  take the maximal
$j$ such that
$(u^j f) v\neq 0$. We have to show that j=0. 
In the contrary case, $(u^j f )v$ is a vector of norm 0:
$((u^jf)v,(u^jf)v)= \pm (v,( u^je)(u^jf)v)=  \pm (v,(u^{2j}h)v)=0$, since
otherwise $(u^{2j}h)v\neq 0$, hence \linebreak 
$(u^{2j}f)v   \neq   0$  (by applying $e$ to it). Hence we
got a non-zero vector of zero norm, unless the module is actually a
$\widehat{gl}(\infty, \Bbb C)$-module. 
\qed\enddemo

Therefore, using Lemma 6.10, Lemma 6.2 and Corollary 4.9,  we have:

\proclaim{Lemma 6.11} If $L(\lambda)$ is a unitary quasifinite
$W_{\infty}$-module, then $\phi_{\lambda}(x) = $\linebreak $\sum_i n_i
e^{r_i x}$, with  $n_i\in \Bbb C$, $r_i\in \Bbb R$.
\endproclaim

Now we can formulate the main result of this section, that follows (in the
same way as Theorem 5.2 in [KR1]), from  Theorem 6.7, Propositions 6.4 and
6.6, and Lemma 6.11:

\proclaim{Theorem 6.12} (a) Let $L(\lambda)$ be  a non-trivial  quasifinite
$W_{\infty}$-module. For each $0\leq \alpha <1$, let $E_{\alpha}$ denote
the set of exponents of $L(\lambda)$ that are congruent to $\alpha$ {\rm
mod }$\Bbb Z$. 
 Then  $L(\lambda)$ is unitary if and only if the following three
conditions are satisfied:

1. All exponents are real numbers.

2. The multiplicities of the exponents $r_i\in E_0$ are positive integers.

3. For exponents $r_i\in E_{\alpha}$ ($0<\alpha <1$), all multiplicities
$mult(r_i)$ are integers and only one of them is negative, $m_{\alpha}:= -
\sum_{r_i\in E_{\alpha}} (mult(r_i))$ is a positive integer, and $r_i -r_j
\leq m_{\alpha}$ for all $r_i, r_j\in E_{\alpha}$.

(b) Any  unitary quasifinite $W_{\infty}$-module $L(\lambda)$  is obtained
by taking tensor product of  unitary irreducible quasifinite highest weight
modules over $\frak g _{s_i}$, $i=1,\dots , m$,  and restricting to
$W_{\infty}$ via the embedding $\hat{\varphi}_{\vec s}$, where $\vec
s=(s_1,\dots,s_m)\in \Bbb R^m$, $0\leq s_i <1$  and $s_i-s_j\notin\Bbb Z$
if $i\neq j$.

\endproclaim

\vskip .4cm

\noindent{\it Acknowledgment.} This research was supported in part by NSF
grant DMS-9622870,  Consejo Nacional de Investigaciones Cient\'\i ficas y
T\'ecnicas, and Secretr\'\i a de Ciencia y T\'ecnica  (Argentina). J.
Liberati would like to thank C. Boyallian for constant help
and encouragement throughout the development of this work, and MIT for its
hospitality.

\vskip .4cm

\

\enddocument